\input amstex
\documentstyle{amsppt}
\magnification=1200

\vsize19.5cm
\hsize13.5cm
\TagsOnRight
\pageno=1
\baselineskip=15.0pt
\parskip=3pt

\def\phi{\varphi}

\nologo
\NoRunningHeads

\topmatter

\title { Perelman's W-functional and stability of  K\"ahler-Ricci flow }\endtitle


\author {Gang $\text{Tian}^{*}$ \& Xiaohua $\text{Zhu}^{**}$}\endauthor
\thanks {
* Partially supported by NSFC  11331001, ** by   NSFC 11331001 and  NSFC 11771019}
 \endthanks \subjclass Primary: 53C25; Secondary: 32J15, 53C55,
  58J05\endsubjclass

\address{\flushpar  *BICMR and School of Mathematical Sciences, Peking University,
Beijing, 100871, China
\newline ** School of Mathematical Sciences, Peking University,
Beijing, 100871, China}\endaddress

\email{\flushpar tian\@math.princeton.edu
\newline  xhzhu\@math.pku.edu.cn} \endemail

\endkeywords

\abstract
 In this expository note, we study the second variation of Perelman¡¯s entropy on the space of Kahler metrics at a K\"ahler-Ricci soliton. We prove that the entropy is stable in the sense of variations. In particular, Perelman¡¯s entropy is stable along the K\"ahler-Ricci flow. The Chinese version of this note has appeared in a volume in honor of professor K.C.Chang (Scientia Sinica Math., 46 (2016), 685-696).

\endabstract
\endtopmatter

\document
\subheading{0. Introduction}

In this note we study the second variation of Perelman's
entropy on the space of K\"ahler metrics  at a K\"ahler-Ricci
soliton. A K\"ahler metric $g_{KS}$ on a
compact manifold $M$ is called a (shrinking) K\"ahler-Ricci
soliton if its K\"ahler form $\omega_{KS}$ satisfies the equation
$$\text {Ric}(\omega_{KS})-\omega_{KS}=L_X\omega_{KS},$$
where $\text {Ric}(\omega_{KS})$ is the Ricci form of $g_{KS}$ and
$L_X\omega_{KS}$ denotes the Lie derivative of $\omega_{KS}$ along
a holomorphic vector field $X$ on $M$. If $X=0$, then $g_{KS}$ is
a K\"ahler-Einstein metric with positive scalar curvature.    We
will  show that the second variation of Perelman's entropy  is
non-positive  in  the space of K\"ahler metrics with $2\pi c_1(M)$
as  K\"ahler class.  Furthermore, if  $(M,g_{KS})$ is a
K\"ahler-Einstein manifold, then  the second variation is
non-positive in the space of K\"ahler metrics with K\"ahler
classes cohomologous to  $2\pi c_1(M)$ ( complex structures on $M$
may vary).  This implies that Perelman's W-functional is stable in
the sense of variations. We will also  determine  the kernel of elliptic
operators which arise in the second variation.   The result of this problem was discussed in the lectures by the
first named author in the Clay summer school on Ricci Flow and
Geometrization  in the summer of 2005.

The organization of this paper is as follows: In Section 1, we
review Perelman's W-functional  and give  a formula of Perelaman's entropy   for the  second
variation. In Section 2, we compute the second variation of
Perelman's entropy  on the space of K\"ahler metrics with  K\"ahler class
$2\pi c_1(M)$ on a fixed complex manifold $M$.  In Section 3, we
extend our calculations to possible  varying  complex structures
on a K\"ahler-Einstein manifold.

\subheading {1. The second variation formula of $\lambda(g)$}

 Recently,  G. Perelman  introduced a functional on a compact
differential manifold $M$ of dimension $n$ [Pe],
$$W(g,f,\tau)=(4\pi\tau)^{-n}\int_M[\tau(R(g)+
|D f|^2)+f-n] e^{-f}dV_g,\tag 1.1$$
 where $R(g)$ denotes the scalar curvature of a Riemannian metric  $g$, f
 is a smooth function $f$ and   $\tau$ is  a
constant. Furthermore,   we may normalize the triple $(g,f,\tau)$
so that
 $$(4\pi\tau)^{-\frac{n}{2}}\int_M e^{-f} dV_g\equiv 1.$$
In our case, we will further normalize the volume of $g$, i.e.,
$$   (2\pi)^{-\frac{n}{2}} \int_M dV_g =   \equiv 1.\tag 1.2$$
Then the $W$-functional can  be reduced to the following
functional  on a pair $(g,f)$,
$$W(g,f)=\int_M[(R(g)+
|D f|^2)+f] e^{-f}dV_g,\tag 1.3$$
 where $(g,f)$ satisfies
 $$   (2\pi)^{-\frac{n}{2}}     \int_M e^{-f} dV_g=    (2\pi)^{-\frac{n}{2}} \int_M dV_g\equiv 1.\tag 1.4$$

For any Riemannian metric $g$  with normalized volume (1.2), we
define Perelman's  entropy as follows,
$$\lambda(g)=\inf_f \{W(g,f)|~ f~\text{satisfies}~(1.4)\}.$$
The number $\lambda(g)$ can be attained by some $f$ (cf. \cite{Ro}). In fact, such
a $f$ is a solution of the equation,
$$2\triangle f+f-|D f|^2+R=\lambda(g).\tag 1.5$$

As in [Pe], we have  the first variation of $\lambda(g)$,
$$\delta\lambda(g)=-   (2\pi)^{-\frac{n}{2}}    \int_M<\delta g, \text{Ric}(g)-g+D^2f>e^{-f}
dV_g,\tag 1.6$$
 where $\text{Ric}(g)$ denotes the Ricci tensor of $g$ and
 $D^2f$ is the  Hessian  of $f$.  It follows from (1.5)  that
  $g$ is a critical metric of $\lambda(g)$ if and only if $g$ is a gradient (shrinking) Ricci
  soliton.  Namely, the metric $g$ satisfies,
   $$\text{Ric}(g)-g=-D^2f,$$
for some smooth function $f$.  Usually, $f$ is called the defining function of Ricci
  soliton $g$.  In fact, one can prove that  $f$ is a unique solution of the equation (1.5)  modulo a constant if  the Riemannian metric
  $g$ is a gradient  Ricci soliton (cf. \cite{TZ4}).

By the standard computation, one can easily get   the second
variation of $\lambda(g)$ at a critical point, i.e.,  a gradient
shrinking Ricci soliton $g_{RS}$. This is given in the following
proposition (also see \cite{CHI, DWW}).

\proclaim {Proposition 1.1} Let $(g_{RS},f)$ be a  gradient
shrinking Ricci soliton on $M$. Let  $L_0$  and $L'$  be  defined
on the space of symmetric tensors of rank 2, respectively, by
$$L_0h=-\frac{1}{2}D^*D h +\text{Rm}(h,\cdot)+\frac{1}{2}(D^2f\cdot h+ h\cdot D^2f)$$
and
$$ \aligned L'h=&\triangle(\text{tr}(h))+\text{tr}(h)+<h,D^2f>_{g_{RS}}-\text{div}(\cdot\text{div}h)\\
 &+ <\text{div} h, D f>_{g_{RS}}-\frac{1}{2}<D (\triangle(\text{tr}(h))), D f>_{g_{RS}}
 -<Df,Df>_{h},\endaligned\tag 1.7$$
  where $<Df,Df>_{h}=\sum_i h^{ij}f_if_j$ and $h^{ij}=\sum_{kl}
g^{ik}g^{jl}h_{kl}$. Then we have
$$\delta^2 \lambda(g)(h,h)=(h,Lh)_{g_{RS}}=\int_M
<h,Lh>_{g_{RS}}dV_{g_{RS}},\tag 1.8$$
 where $L$ is  defined by
$$Lh=L_0h+\text{div}^*\cdot\text{div}h+\frac{1}{2}D^2(\text{tr}(h))-
D^2(P^{-1}\cdot L'(h)),$$
  and $P$ is given by
 $$P\psi=2\triangle \psi+ 2<Df, D\psi>_{g_{RS}}+\psi, ~~~\forall ~\psi\in C^{\infty}(M).$$

\endproclaim

\demo{Proof} According to [Be], we have
$$\aligned \delta \text{Ric}(g)(h)=&\frac{1}{2}D^*D h
-\text{Rm}(h,\cdot)+\frac{1}{2}(\text{Ric}\cdot
h+h\cdot\text{Ric})\\
&-\text{div}^*\cdot\text{div}h-\frac{1}{2}D^2(\text{tr}(h)).\endaligned\tag 1.9$$
 On the other hand, by (1.5), one gets an equation for $\delta f$,
 $$P (\delta f)= L'h.$$
Combining these two, we obtain (1.8).\qed
\enddemo

\subheading {2. The   case  for a fixed  complex structure}

In this section, we compute the second variation of $\lambda(g)$
restricted to the space of  K\"ahler metrics  with K\"ahler forms
in  $2\pi c_1(M)>0$. This variation is computed at  an $n$-dimensional  (shrinking)
K\"ahler-Ricci soliton and contains more information than that in
the real case.  In  K\"ahler case, $W$-functional can be rewritten as
$$W(g,f,\tau)=\int_M[\tau(R(g)+
|\bar\partial f|^2)+f] e^{-f}dV_g,\tag 2.1$$
where
$$\int_M e^{-f} \omega_g^n =  \int_M  \omega_g^n   =(2\pi)^{n} c_1(M)^n.$$

 Let $(g_{KS},X)$ be  a K\"ahler Ricci
soliton with its K\"ahler form $\omega_{KS}$ in  $2\pi c_1(M)$ on
a compact complex manifold  $(M,J)$, where $J$ denotes a complex
structure of $M$.
We consider all  K\"ahler metrics $g$ with
K\"ahler forms in $2\pi c_1(M)$.
Without loss of generality, we
may assume that $g=g_t=g_{KS}+th_0$ is a family of such K\"ahler
metrics, where
  $$h_0=\sum_{ij}
\partial_i\partial_{\overline j}\psi dz^i\otimes d\overline
z^j$$
 for some real-valued smooth function $\psi$.
 We shall compute $\frac{d^2\lambda(g_t)}{dt^2}|_{t=0}$.
Put
 $$P_0\psi=2\triangle\psi+\psi-(X+\overline X)(\psi),$$
$$L_1\psi=\triangle\psi+\psi-X(\psi),$$
and
$$L_1'\psi =\triangle\psi-X(\psi).$$

 \proclaim{ Lemma 2.1} Let $g_t=g_{KS}+th_0$ be a family of
 K\"ahler metrics as above and $f=f_t$ be a family of smooth
functions which are solutions of  (1.5) associated to  $g_t$. Let
$u=\frac{d f}{dt}|_{t=0}$. Then
$$ P_0(u-X(\psi))= (L_1'\cdot
L_1)(\psi).$$
\endproclaim

\demo{Proof} Differentiating  relations (1.5) at $t=0$, we have
$$\aligned  P_0(u)&=2\triangle u+u- 2\text{re}(X(u))\\
&=\triangle^2
\psi+<\text{Ric}(g_{KS}),\sqrt{-1}\partial\partial\psi>+\psi_{i\overline
j}(2f_{j\overline i}-f_jf_{\overline i}).\endaligned  \tag 2.2$$
 It follows
$$\aligned
P_0(u)&=\triangle^2 \psi+\triangle\psi+\psi_{i\overline
j}(f_{j\overline i}-f_jf_{\overline i})\\
&=\triangle^2
\psi+\triangle\psi+\triangle(X(\psi))-X(\triangle\psi)-X(\overline{X(\psi)})\\
&=(\triangle-X)[\triangle
\psi+\psi-X(\psi)] +2\triangle(X(\psi))+X(\psi)- X(X(\psi))-\overline X(X(\psi))\\
&=[(\triangle-X)\cdot L_1](\psi)+P_0(X(\psi)).\endaligned$$
 The lemma follows.\qed\enddemo

\proclaim{Proposition 2.1} Let $g_t=g_{KS}+th_0$ be a family of
K\"ahler metrics as in Lemma 2.1. Then
$$\frac{d^2\lambda(g_t)}{dt^2}|_{t=0}
= \int_M \psi \times [P_0^{-1}\cdot (\overline
L_1'L_1')\cdot(\overline L_1L_1)](\psi) e^{-f}\omega_{KS}^n\le
0.\tag 2.3$$
 Moreover the
equality in (2.1) holds if and only if
$\psi=\theta_v+\overline\theta_v$ for some holomorphic vector
field $v$ on $M$, where
 $\theta_v$ is a potential  associated to $v$ defined by
$$i_v(\omega_{KS})=\sqrt{-1}\overline\partial\theta_v.$$
\endproclaim

\demo {Proof} First we see
 $$\frac{d\lambda(g_t)}{dt}=-\int_M<\sqrt{-1}\partial\overline\partial\psi+\theta,
 \text{Ric}(\omega)-\omega+\sqrt{-1}\partial\overline\partial f>e^{-f}\omega^n.$$
 Since
 $$\frac{d\text{Ric}(g_t)}{dt}|_{t=0}=-\sqrt {-1}\partial\overline\partial(\triangle\psi),$$
 we have
$$\aligned &\frac{d^2\lambda(g_t)}{dt^2}|_{t=0}\\&=\int_M<\partial\overline\partial\psi,
\partial\overline\partial(\triangle\psi+\psi-\frac{df_t}{dt}|_{t=0})>
e^{-f}\omega_{KS}^n \\
&=\int_M \sum \psi_{i\overline
j}(\triangle\psi+\psi-u)_{j\overline
i}e^{-f}\omega_{KS}^n.\endaligned\tag 3.4$$
  Integrating  by parts, we get
$$\aligned &\frac{
d^2\lambda(g_t)}{dt^2}|_{t=0}\\&= \int_M
(\triangle\psi+\psi-u)[\triangle^2\psi+X(\overline{X(\psi)})-2\text{re}(X(\triangle\psi))-\psi_{i\overline
j}f_{j\overline i}] e^{-f}\omega_{KS}^n\\
&=\int_M (\triangle\psi+\psi-u)[\triangle(\triangle\psi- X(\psi))
-\overline X(\triangle\psi-X(\psi))]
e^{-f}\omega_{KS}^n\\
&=\int_M (\triangle\psi+\psi-u)[\overline {(\triangle-X)
(\triangle\psi-\overline{X(\psi)})}]
e^{-f}\omega_{KS}^n\\
&= \int_M
(\triangle-X)(\triangle\psi+\psi-u)\times(\triangle\psi-X(\psi))e^{-f}\omega_{KS}^n\\
&= \int_M[
(\triangle-X)(\triangle\psi+\psi-X(\psi))+(\triangle-X)(X(\psi)-u)]
\\
&\times(\triangle\psi-X(\psi))e^{-f}\omega_{KS}^n.
\endaligned\tag 2.5$$
 By Lemma 2.1, we derive from (2.5),
$$\aligned &\frac{
d^2\lambda(g_t)}{dt^2}|_{t=0}\\&=\int_M[P(u-X(\psi))+(\triangle-X)(X(\psi)-u)]
\times(\triangle\psi-X(\psi))e^{-f}\omega_{KS}^n\\
&=\int_M [(\triangle-\overline X )
(u-X(\psi))+ (u-X(\psi))]\times(\triangle\psi-X(\psi))e^{-f}\omega_{KS}^n\\
&=\int_M \overline L_1(u-X(\psi))\times(\triangle\psi-X(\psi))e^{-f}\omega_{KS}^n\\
&=\int_M  (u-X(\psi))\times [L_1\times (\triangle- X )](\psi) e^{-f}\omega_{KS}^n\\
& =  \int_M P_0^{-1}((L_1'L_1)\psi) \times (L_1'L_1)(\psi)
e^{-f}\omega_{KS}^n\\
& =  \int_M [\overline{(L_1'L_1)}\cdot P_0^{-1}\cdot
((L_1'L_1)](\psi) \times \psi e^{-f}\omega_{KS}^n.\endaligned$$
 Since  any two operators of $P_0,L_1,L_1'$  commute,   $P_0^{-1}$
  commutes with $L_1,L_1'$. Thus we have
$$\frac{
d^2\lambda(g_t)}{dt^2}|_{t=0}= \int_M \psi \times [P_0^{-1}\cdot
(\overline L_1'L_1')\cdot(\overline
L_1L_1)](\psi)e^{-f}\omega_{KS}^n.$$
 Note that  $P_0, L',\overline L'$ are all elliptic, so does
  $P_0^{-1}\cdot (\overline L_1'L_1')$. This shows that
  $$\frac{
d^2\lambda(g_t)}{dt^2}|_{t=0}\le 0$$
 and the equality holds if and only if
 $$\overline L_1L_1(\psi)=0.$$
   Then  Proposition 2.1 will be  completed from the next lemma. \qed
\enddemo

\proclaim{Lemma 2.2} The operator $\overline L_1L_1$ is  real and
nonnegatively definite. Moreover, there is an isomorphism  between
$\text{ker}(\overline L_1L_1)$ and  the linear  space  $\eta(M)$
of holomorphic vector fields on $(M,J)$ given by a relation
$\psi=\theta_v +\overline\theta_v$ for some $v\in\eta(M)$.
\endproclaim

\demo{Proof} It suffices  to prove the second part of the Lemma.
This follows from an argument   in Appendix of [TZ1]. In  fact
$$\overline L_1L_1\psi=0$$
 implies that
 $$L_1\psi= \theta_u$$
 for some $u\in\eta(M)$.  From the proof of Lemma A.2 in [TZ1], we see that
 $$\psi=\theta_v+\overline{\theta_{v'}}$$
  for some $v,v'\in \eta(M)$. Since $\psi$ is  a real-valued
  function,  $v'$ must be equal to $v$. Thus the lemma is true.
  \qed\enddemo

The formula (2.3) in Proposition 2.1 can be generalized to the
variation of  any K\"ahler metrics $g$  for the fixed complex
structure if the underlying manifold $M$ is  K\"ahler-Einstein
with positive scalar curvature as follows. Let $g_{KE}$ be a
K\"ahler-Einstein metric on $(M,J)$ and
$$g_t=g_{KE}+t(\theta+\sum_{ij}
\partial_i\partial_{\overline j}\psi dz^i\otimes d\overline
z^j)\tag 2.6$$
  be a family of  K\"ahler metrics with
  $$\int_M \omega_{g_t}^n=\int_M \omega_{KE}^n,\tag 2.7$$
 where $\theta$ be a hermitian and symmetric tensor with respect  to the
complex structure $J$. It is easy to see that  condition (2.7)
means
$$\int_M \text{tr}_{\omega_{KE}}(\theta)\omega_{KE}^n=0.$$
Without loss of generality, we may further assume that the
corresponding (1,1)-form of $\theta$ is harmonic associated to the
metric $\omega_{KE}$. This implies that
$$ d [\text{tr}_{\omega_{KE}}(\theta)]=0.$$
Thus we get
  $$\text{tr}_{\omega_{KE}}(\theta)=0.\tag 2.8$$

 \proclaim{Proposition 2.2} Let $(M,J)$ be  a K\"ahler-Einstein manifold with positive first Chern class
  and $g_{KE}$ be a K\"ahler-Einstein metric on $M$. Let $g_t$ be
 a family of  K\"ahler metrics of the form (2.4) and satisfying (2.6).
 Then
$$\aligned &\frac{d^2\lambda(g_t)}{dt^2}|_{t=0}\\
&= \int_M \|\theta\|^2\omega_{KE}^n+\int_M(<\Cal D^*\Cal D\psi,
(P_0')^{-1}(\Cal D^*\Cal D\psi)>\omega_{KE}^n,\endaligned\tag
2.9$$
 where  $P_0'$ and  $\Cal D$  are defined, respectively, by
$$P_0'\psi= 2\triangle_{KE}\psi+\psi$$
and
$$\Cal D\psi= \sum\psi_{\overline i\overline j}dz^{\overline
i}dz^{\overline j},$$
 and  $\Cal D^*$ is the adjoint operator of $\Cal D$. \endproclaim

\demo {Proof} By
 $$\frac{d\lambda(g_t)}{dt}=-\int_M<\sqrt{-1}\partial\overline\partial\psi+\theta,
 \text{Ric}(\omega)-\omega+\sqrt{-1}\partial\overline\partial f>e^{-f}\omega^n$$
 and
 $$\frac{d\text{Ric}(g_t)}{dt}|_{t=0}=-\sqrt {-1}\partial\overline\partial(\text{tr}\theta+\triangle\psi)
 =-\sqrt {-1}\partial\overline\partial(\triangle\psi),$$
 we have
$$\aligned &\frac{d^2\lambda(g_t)}{dt^2}|_{t=0}\\&=\int_M<\partial\overline\partial\psi+\theta,
\partial\overline\partial(\triangle\psi+\psi-\frac{df_t}{dt}|_{t=0})+\theta>
\omega_{KE}^n.\endaligned\tag 2.10$$
 Taking  the integral by parts, we get
$$\aligned &\frac{d^2\lambda(g_t)}{dt^2}|_{t=0}=\int_M\|\theta\|^2\omega_{KE}^n\\
&+\int_M<\partial\overline\partial\psi,
\partial\overline\partial(\triangle\psi+\psi-\frac{df_t}{dt}|_{t=0})>\omega_{KE}^n.
\endaligned\tag 2.11$$

Since $f=f_t$ satisfies
$$2\triangle f+f-|D f|^2+R=\lambda(g),$$
 differentiating at $t$ on the both sides, we get
 $$\aligned
P_0'\frac{df_t}{dt}|_{t=0}&=2\triangle\frac{df_t}{dt}|_{t=0}+\frac{df_t}{dt}|_{t=0}\\
&=\triangle^2\psi+\triangle\psi+\triangle(\text{tr}\theta)+\text{tr}\theta\\
 &=\Cal D^*\Cal D\psi.\endaligned\tag 2.12$$
 Note that $P_0'$ is invertible since the first non-zero eigenvalue is $1$
  on the K\"ahler-Einstein manifold.
  Using this relation, we obtain
$$\aligned &\int_M<\partial\overline\partial\psi,
\partial\overline\partial(\triangle\psi+\psi-\frac{df_t}{dt}|_{t=0})>\omega_{KE}^n\\
&=\int_M(<\Cal D^*\Cal D\psi, (P_0')^{-1}(\Cal D^*\Cal
D\psi)>\omega_{KE}^n.\endaligned\tag 2.13$$
 Thus combining (2.11) and (2.13), we prove (2.9).
\qed
\enddemo

The following corollary shows that Proposition  2.1 is not true in
general if K\"ahler metrics are not fixed in the K\"ahler class
$2\pi c_1(M)$.

 \proclaim {Corollary 2.1}Let $(M,J)$
be a K\"ahler-Einstein manifold with positive first Chern class.
Suppose that $\text{dim}H^{1,1}(M,J)\ge 2.$ Then
$\delta^2\lambda(g)(h,h)$ is not non-positive at a
K\"ahler-Einstein metric $g_{KE}$ for the variation of general
K\"ahler metrics and so $g_{KE}$ is not a local maximum of
$\lambda(g)$ in the total space of K\"ahler metrics.
\endproclaim

\demo{Proof} Let $\omega'$ be another harmonic (1,1)-form of
$(M,J)$ which is not a multiple of $c_1(M)$. Then there are two
number $a$ and $b$ such that
$$ an+ b \text{tr}_{\omega_{KE}}\omega'=0.$$
Let $\theta=a\omega'+b\omega_{KE}$ and $h_0$  its corresponding
hermitian and symmetric tensor. Then
$$\int_M
\text{tr}_{\omega_{KE}}(\theta)\omega_{KE}^n=0.$$
  Thus by Proposition 2.2, we have
  $$\delta^2\lambda(g)(h_0,h_0)= \int_M \|\theta\|^2\omega_{KE}^n>0.$$
This implies that $\delta^2\lambda(g)$ is not non-positive in the
direction of $h_0$, so the corollary is true.\qed \enddemo

\subheading {3. The case for varying complex structures}

 In this section, let $(M, g_{KE}, J_0)$ be a K\"ahler-Einstein manifold,  we
 will study  the second variation of $\lambda(g)$ at $g_{KE}$ when restricted   K\"ahler
 metrics with K\"ahler forms  cohomologous to  $2\pi c_1(M)$. Set
$$\aligned \Cal W=\{h|~&\text{there is a family of K\"ahler metrics}~ (g_t,J_t)(0\le t\le\epsilon)
~\text{such that}\\
&h=\frac{dg_t}{dt}|_{t=0},~(g_0,J_0)=(g_{KE},J_0),~\text{and}~
[\omega_t]= 2\pi c_1(M)\}.
\endaligned\tag 3.1$$
Here $J_t$ denotes a family of complex structures on $M$ and
$\omega_t$ denotes the K\"ahler form of $g_t$ . We shall prove

 \proclaim{Theorem 3.1}  The operator $L$ defined in Proposition 1.1 is non-positive on $\Cal W$.
Namely,  for any $h\in \Cal W$, we have
$$ \delta^2 \lambda(g_{KE})(h)=(h,Lh)_{g_{KE}}=\int_M<h,Lh>_{\omega_{KE}}\omega_{KE}^n\le 0.\tag 3.2$$
Moreover, there exists an isomorphism
 $$\imath: \text{ker}(L)\to \eta(M,J_0)+H^1(M,J_0,\Theta),$$
 where $\text{ker}(L)$ denotes the kernel of  $L$, and $\eta(M,J_0)$
  is the space of holomorphic vector fields associated to
 the complex structure $J_0$ on $M$ and  $H^1(M, J_0,\Theta)$
 is the $\overset {\vee} \to C$ech cohomology class associated to the infinitesimal
deformation of complex structures on $M$ [Kod].
\endproclaim

\proclaim{Remark 3.1}  According to Corollary 2.1 in Section 2, we
see that  Theorem 3.1 is not true in general for  K\"ahler metrics
without the assumption that the K\"ahler class $2\pi c_1(M)$  is
fixed.  The same observation was made in [CHI] in the Riemannian
case.
\endproclaim

As before, we let $J_t (0\le t\le\epsilon)$ be a family of complex
structures on a K\"ahler-Einstein manifold $M$ with $J_0=J $. Then
according to [Koi], one can decompose $\frac{dJ_t}{dt}|_{t=0}$ as
a direct product into
$$\frac{dJ_t}{dt}|_{t=0}=L_ZJ + I_E$$
with
$$\int_M <L_ZJ,  I_E>_{\omega_{KE}}\omega_{KE}^n=0,$$
 where $L_Z$ denotes the Lie derivative along a vector field
$Z$ on $M$ and $I_E$ is the part of an essential infinitesimal
deformation of complex structures on $M$. If we let $h'$ be a
covariant tensor of rank 2  defined by
 $$h'(X,Y)=\omega_{KE}(X,I_EY),\tag 3.3$$
  then  $h'$ is anti-hermitian, and so
it is a real part of some  $(0,2)$-type tensor $I=I_{\overline
i\overline j} d\overline z^i\otimes d\overline z^j$, i.e.,
$$ h'= \text{Re} (I).$$
 Moreover,  $I$ satisfies ([Koi]),
$$ \nabla_{\overline k} I_{\overline i\overline j}=\nabla_{\overline j} I_{\overline i\overline
 k},~\forall ~i,j,k,\tag 3.4$$
 and
 $$\sum_j \nabla^{\overline j} I_{\overline i\overline j}=0,~\forall ~i.\tag 3.5$$
The relations (3.4) and (3.5) imply that the complexification of
$I_E$ is a $\overline \partial$-closed, $(0,1)$-form with  values
in $\Theta_{J_0}$.  This defines a  $\overset {\vee} \to C$ech
cohomology class in $H^1(M,\Theta_{J_0})$, where $\Theta_{J_0}$
denotes the $(1,0)$-typed  tangent sheaf  associated to $J_0$ on
$(M,J_0)$ [Kod]. (3.5) also implies that
$$\text {div} h'=\sum_{\alpha} D_{e_{\alpha}}
h'(\cdot, e_{\alpha})=0.\tag 3.6$$

 Let $\rho_t$ be an one-parameter diffeomorphisms group generated by the vector
 $-Z$ and  $\rho_t^* g_t$ be a family of induced Riemannian metrics, where
$$g_t(X,Y)=\omega_t(X,J_tY).$$
Then
$$\aligned (\rho_t)^*g_t(X,Y)&= g_t((\rho_t)_*X,(\rho_t)_*Y)=
\omega_t((\rho_t)_*X,J_t((\rho_t)_*Y))\\
&=(\rho_t)^*\omega_t(X,((\rho_t)^*J_t)Y).\endaligned$$
  It follows
 $$\aligned \tilde h(X,Y)&=\frac{d[(\rho_t)^*g_t]}{dt}|_{t=0}(X,Y)\\
 &=\frac{d
 [(\rho_t)^*\omega_t]}{dt}|_{t=0}(X,JY)+ \omega_{KE}(X,\frac{d[(\rho_t)^*J_t]}{dt}|_{t=0}Y) \\
&=\frac{d [(\rho_t)^*\omega_t]}{dt}|_{t=0}(X,JY)
+\omega_{KE}(X,I_EY).\endaligned $$

 Set
$$\aligned \Cal W_0=&\{h~\text{is a covariant symmetric tensor of
rank 2  such that}\\ & Lh=0~\text{and} ~\text{div}(
h)=0\}.\endaligned$$

\proclaim{Lemma 3.1} Let $h'$ be the covariant 2-tensor in (3.3).
Then $h'\in \Cal W_0$.
\endproclaim

\demo {Proof}  By a direct computation, it was showed in [Koi]
that (3.4) and (3.5) implies,
 $$ L_0 h'=\frac{1}{2}D^*D h' -\text{Rm}(h',\cdot )=0. \tag 3.7 $$
On the other hand, one can decompose $h'$ into a symmetric part
$b$ and an anti-symmetric part $a$ which is orthogonal in the
sense of inner product
 $$(a,b)_{\omega_{KE}}=\int_M<a,b>_{\omega_{KE}}\omega_{KS}^n.$$
Since the operator $L_0$ keeps the symmetry  and anti-symmetry,
the anti-symmetric part $a$ of $h'$ also satisfies  equation
(3.7), and consequently, $a$ is parallel, i.e.,
$$Da=0.$$
 By using  the Ricci  identity
 $$ D^iD_i a-D^{\overline i}D_{\overline i} a= 2a,$$
we see that $a=0$.  This implies that  $h'$ is a symmetric
2-tensor. By using  (3.5) and the fact
$\text{tr}_{\omega_{KE}}h'=0$, we also get
$$ Lh'=L_0 h'=0.$$
Hence, we have $h'\in W_0$. \qed\enddemo

 \proclaim{Lemma 3.2}  Assume that $\omega_t\in 2\pi c_1(M)$.
 Then there is a  smooth real-valued function $\psi$ on $M$ such that
 $$\frac{d[(\rho_t)^*\omega_t]}{dt}|_{t=0}=\sqrt{-1}\partial\overline\partial\psi.$$
\endproclaim

\demo{Proof} Let  $\theta_1$ and $\theta_2$ be  in $A^{1,1}(M,J)$
and $A^{2,0}(M,J)$ respectively, such that
$$\frac{d[(\rho_t)^*\omega_t]}{dt}|_{t=0}=\theta_1+\text{Re}(\theta_2).$$
Clearly,  tensor $h_1$ defined by $h_1(X,Y)=\theta_1(X,JY)$ is
symmetric and hermitian and tensor $h_2$ defined by
$h_2(X,Y)=\text{Re}(\theta_2)(X,JY)$ is anti-symmetric. Since $h$
and $h'$ are both symmetric according to $h\in \Cal W$ and Lemma
3.1, we see that $\text{Re}(\theta_2)$ must vanish. Thus we get
 $$\frac{d[(\rho_t)^*\omega_t]}{dt}|_{t=0}=\theta_1.$$
Note that $\theta_1$ is also an exact form because of $\omega_t\in
2\pi c_1(M)$. Therefore the lemma is true. \qed\enddemo

According to Lemma 3.1 and 3.2, we see that $h_1$ and $h'$ are
hermitian and anti-hermitian symmetric tensors respectively. Then
$<h_1,h'>_{\omega_{KE}}=0.$ Thus we have

 \proclaim{Lemma 3.3}
$$(h_1,h')_{\omega_{KE}}=0.$$
\endproclaim

 Combining Lemma 3.2 and 3.3, we   get
\proclaim {Proposition 3.1}
$$\Cal W/ \text{diff(M)}\cong A^{1,1}(M,J_0)\bigoplus H^1(M,\Theta_{J_0}).$$
\endproclaim

\demo {Proof of Theorem 3.1}
  By Lemma 3.2, we have
 $$ \tilde h(X,Y)=h_1(X,Y)+h'(X,Y).$$
  Since  $(h_1,h')_{\omega_{KE}}=0$
 by Lemma 3.3,  we see that
  $$(\tilde h,L\tilde h)=(h_1,Lh_1)+(h',Lh')=(h_1,Lh_1).\tag 3.8$$
The last equality follows from Lemma 3.1.  On the other hand, by
Lemma 3.2, we see that there is a smooth real-valued function
$\psi$ such that
$$h_1(X,Y)=\sqrt{-1}\sum_{i,j}\partial_i\partial_{\overline j}\psi
dz^i\wedge d\overline z^j(X,JY).$$
  Then according to Proposition 2.2, we have
$$(h_1,Lh_1)=\frac{d^2\lambda(g_t)}{dt^2}|_{t=0}\le 0,$$
where $g_t$ is a family K\"ahler metrics defined as in  Lemma 2.1
with $h_0$ replaced by $h_1$. Thus
$$(\tilde h,L\tilde h)=(h_1,Lh_1)\le 0.$$
 Since $W$-functional is invariant under
 diffeomorphisms, we obtain
 $$ (h,Lh)=(\tilde h, L\tilde h)\le 0.$$

By relation (3.8) and Proposition 2.1, $\tilde h$ is a  kernel of
L iff $\tilde h=h_1+h',$
 where $h'$ is defined by (3.3) and
 $$h_1=\text{Re}(\sum_{ij} \partial_i\partial_{\overline
j}\psi dz^i\otimes d\overline z^j)$$
 for some real-valued function $\psi$ which satisfies
  $$\psi=\theta_v +\overline\theta_v,$$
where $v\in\eta(M)$.  Thus the operator $L$ induces an injective
homomorphism
 $$\imath: \text{Ker}(L)\to \eta(M,J_0)+H^1(M,\Theta_{J_0}).$$
  It is clear that $\imath$ is surjective by Lemma 3.1 and (3.8) and the fact that $\lambda(g)$ is invariant under
   the holomorphic  transformations. Therefore the
theorem  is true.\qed\enddemo

\proclaim{Remark 3.2} The relation (3.2) can be also obtained by
showing that $g_{KE}$ is a global minimizer of $\lambda(g)$ in the
total space of K\"ahler metrics $(g,J)$ with K\"ahler classes
cohomologous to $2\pi c_1(M)$.  In fact, we have
$$\aligned \lambda(g_{KE})&=W(g_{KE},1)=\int_M \omega_{KE}^n \\
&=W(g,1)\ge\inf_f \{W(g,f)|~
f~\text{satisfies}~(1.4)\}=\lambda(g).\endaligned\tag 3.9$$
  But Theorem 3.1  determines  more explicitly  the kernel of the elliptic operator
which arises in the second variation. In particular,  the dimension of kernel  is finite modula diffeomorphisms group. The latter will be very useful, for example, Sun and Wang used it to study the stability of K\"ahler-Ricci flow \cite{SW}.  We conjecture
that there exists an analogous version of Theorem 3.1 for complex manifolds which admit
K\"ahler-Ricci solitons.

\endproclaim

\vskip10mm

\noindent{\bf Notes.}  The  W-functional was introduced by Perelman in 2012.  This functional plays a crucial role  in his proofs of  Poincar\'e conjecture and  Thurston's Geometrization theorem (cf. \cite{P1, P2, P3,  MT1, MT2}). Our article was  posted in arXiv:0801.3504, 2008 \cite{TZ3}.
After that, there are other  new developments  in the study of  Perelman's W-functional and entropy  with  their applications in Ricci flow, cf.
\cite{W, HT1, HT2, SW, Pa1, Pa2, CW, Ba, TZh}, etc..  For example,   variant  generalizations  of Theorem 3.1.   were  studied for a complex manifold which admits a  K\"ahler-Ricci soliton by other people, cf. \cite {W, HT1, HT2, Pa1, Pa2}. The authors also used   Perelman's  $W$-functional and entropy  to  study  the stability and convergence  of  K\"ahler-Ricci flow  in later papers,
 \cite {Zh, TZ4, TZ5}.

\Refs\widestnumber\key{[Chw]}


\item {[Bam]} Bamler, R., Convergence of Ricci flows with bounded scalar curvature, arXiv:1603.05235.

 \item {[Be]} Besse, A.L., Einstein manifold,
Springer-Verlog, 1987.

\item{[Ca]}Cao, H.D., Deformation of K\"ahler metrics to
K\"ahler-Einstein metrics on compact K\"ahler manifolds, Invent.
Math.,  81 (1985), 359-372.

 \item {[CHI]} Cao, H.D.,  Hamilton, S. and  Ilmanen, T.,
Gaussian densities and stability for some Ricci solitons,
preprint.

\item {[CT1]}  Chen, X.X. and Tian, G., Ricci flow on
        K\"ahler-Einstein  surfaces, Invent. Math.,
        147 (2002), 487-544.

\item {[CT2]}  Chen, X.X. and Tian, G., Ricci flow on
        K\"ahler-Einstein  manifolds,  Duke Math. J., 131 (2006),
        17-73.

\item{[CW]} Chen, X. and Wang, B.,
Space of Ricci flows (\uppercase\expandafter{\romannumeral 2}), arXiv:1405.6797.


\item{[DWW]} Dai, X., Wang, X. and Wei, G., On the stability of
K\"ahler-Einstein metrics, Arxiv Math/0504527

\item{[GW]} Greene, R., Wu, H., Lipschitz convergence of
Riemannian manifolds, Pacific J. Math., 131 (1988), 119-141.

\item {[Ha]}  Hamilton, R.S., Three manifolds with positive Ricci
Curvature, J. Diff. Geom., 17 (1982), 255-306.

\item {[HM1]}  Hall S and  Murphy T.,  On the linear stability of K¡§ahler-Ricci solitons,  Proc Amer Math Soc, 139 (2011),  3327-3337.

\item {[HM2]} Hall S and Murphy T.,  Variation of complex structures and the stability of K¡§ahler-Ricci solitons. Pacific J Math,
265 (2013), 441-454.

\item {[Kod]} Kodaira, K., Complex manifolds and
deformation of complex structures, Springer-Verlag, 1986.

 \item{[Koi]} Koiso, N., Einstein metrics and complex
structures, Invent. Math., 73 (1983), 71-106.

\item{[MT1]}  Morgan J and  Tian G.,  Ricci Flow and the Poincar¡äe Conjecture. Clay Mathematics Monographs, 3. Providence: Amer
Math Soc; Cambridge: Clay Mathematics Institute, 2007.

\item{[MT2]} Morgan J, Tian G.,  The Geometrization Conjecture. Clay Mathematics Monographs, 5. Providence: Amer Math Soc;
Cambridge: Clay Mathematics Institute, 2014.

\item{[Pa1]} Pali, N., Variational stability of K?hler-Ricci solitons,  Adv. Math.  290 (2016), 15-35.

\item{[Pa2]} Pali, N.,  Concavity of Perelman's W -functional over the space of K?hler potentials. Eur. J. Math.  3  (2017),  no. 3, 587-602.

\item{[P1]} Perelman, G., The entropy formula for the Ricci flow
and its geometric applications, ArXiv:0211159, 2002.

\item{[P2]}  Perelman G., Ricci flow with surgery  on  three-manifolds. ArXiv:0303109, 2003.

\item{[P3]}  Perelman G.,  Finite extinction time for the solutions to the Ricci flow on certain three-manifolds. ArXiv:0307245, 2003.

\item{[Ro]} Rothaus, O., Logarithmic Sobolev inequality and the spectrum of Schr¡§odinger operators. J Funct Anal, 42 (1981),
110-120.

 \item{[S1]} Sesum, N., Convergence of a K\"ahler-Ricci
flow, Math. Res. Lett., 12 (2005), 623-632.

\item{[S2]} Sesum, N., Linear and Dynamical stability of Ricci
flat metrics, Duke Math. J., 133 (2006), 1-26.


\item{[SW]} Sun, S. and Wang Y.,  On the K¡§ahler-Ricci flow near a K¡§ahler-Einstein metric. J Reine Angew Math, 699 (2015), 143-158.


\item{[TZ1]} Tian, G. and  Zhu, X., Uniqueness of K\"ahler-Ricci
solitons, Acta Math., 184 (2000), 271-305.

\item{[TZ2]} Tian, G. and  Zhu, X., Convergence of the
K\"ahler-Ricci flow, J. Amer Math. Sci., 17 (2007), 675-699.

\item{[TZ3]} Tian, G. and  Zhu, X.,  Perelman's $W$-functional and stability of K\"ahler-Ricci flows,  arXiv:0801.3504, 2008.

\item{[TZ4]}  Tian, G. and  Zhu X.,  Convergence of the K¡§ahler-Ricci flow on Fano manifolds. J Reine Angew Math, 678 (2013), 223-245.

\item{[TZ5]}  Tian, G., Zhang, S, and  Zhang, Z., Perelman¡¯s entropy and K¡§ahler-Ricci flow on Fano manifolds. Tran Amer Math Soc,
365 (2013),  6669-6695

\item{[TZh]} Tian, G. and Zhang, Z.L.,
Regularity of K\"ahler-Ricci flows on Fano manifolds, Acta Math., 216 (2016), 127-176.

\item{[W]} Wang, Y.,  Second variation of -functional at K¡§ahler-Ricci solitons: On Ricci solitons and Ricci flows. PhD Thesis.
Madison: University of Wisconsin, 2011.

\item{[Zh]} Zhu X.,  Stability of K¡§ahler-Ricci flow on a Fano manifold. Math Ann,  356 (2013),  1425-1454.

\endRefs
\enddocument